\documentclass[12pt]{amsart}
\usepackage{verbatim}

\newtheorem{thm}{Theorem}
\newtheorem{prop}{Proposition}
\newtheorem{conj}{Conjecture}
\newtheorem{lemma}{Lemma}
\newtheorem{cor}{Corollary}
\theoremstyle{definition}
\newtheorem{defn}{Definition}
\newtheorem{rk}{Remark}
\newtheorem{rks}[rk]{Remarks}

\newcommand{\C}{\mathbf{C}}
\newcommand{\F}{\mathbf{F}}
\renewcommand{\L}{\mathbb{L}}
\renewcommand{\P}{\mathbb{P}}
\newcommand{\Q}{\mathbf{Q}}

\newcommand{\Z}{\mathbf{Z}}

\newcommand{\sM}{\mathcal{M}}
\newcommand{\sN}{\mathcal{N}}

\newcommand{\iso}{\overset{\sim}{\to}}
\newcommand{\1}{\mathbf{1}}

\newcommand{\car}{\operatorname{char}}
\newcommand{\kim}{{\operatorname{kim}}}
\newcommand{\num}{{\operatorname{num}}}
\newcommand{\rat}{{\operatorname{rat}}}
\newcommand{\alg}{{\operatorname{alg}}}
\newcommand{\End}{\operatorname{End}}
\newcommand{\Hom}{\operatorname{Hom}}
\newcommand{\Pic}{\operatorname{Pic}}
\numberwithin{equation}{section}

\begin{document}
\title{Motivic zeta functions of motives}
\author{Bruno Kahn}
\address{Institut de Math\'ematiques de Jussieu\\175--179 rue du Chevaleret\\75013
Paris\\France}
\email{kahn@math.jussieu.fr}
\date{June 16, 2006}
\maketitle

\section*{Introduction} Let $\sM$ be a tensor category with coefficients in a field $K$ of
characteristic $0$, that is, a $K$-linear pseudo-abelian symmetric monoidal category such that
the tensor product $\otimes$ of $\sM$ is bilinear. Then symmetric and exterior powers of an
object $M\in \sM$ make sense, by using appropriate projectors relative to the action of the
symmetric groups on tensor powers of $M$. One may therefore introduce the \emph{zeta function}
of $M$, a power series with coefficients in
$K_0(\sM)$
\cite[\S 13.3.2]{andre}:
\[Z_\sM(M,T) = \sum_{n\ge 0} [S^n(M)]T^n.
\]

Yves Andr\'e showed that this zeta function is rational when $\sM=\sM(k)$ is the category of
pure motives for some adequate equivalence relation over a field $k$ and $M$ is
finite-dimensional in the sense of Kimura-O'Sullivan (ibid., Prop. 13.3.3.1). He raised the
question of functional equations, in the light of Kapranov's result for curves \cite{kapranov};
this was achieved by Fran\-zis\-ka Heinloth in \cite{heinloth} when $M$ is the motive of an
abelian variety. Heinloth also proved in
\cite[Lemma 4.1]{heinloth} that Andr\'e's rationality argument in fact applies to any
finite-dimensional object in any $K$-linear tensor category $\sM$.

In the same vein, we propose in Theorem \ref{t1} below a functional equation for any
finite-dimensional
$M\in\sM$ in any rigid $K$-linear tensor category; this theorem has a very simple proof. In
particular, modulo homological equivalence we get a functional equation for the motivic zeta
function of any smooth projective variety $X$ for which the sign conjecture holds, that is, the
sum of the even K\"unneth projectors of $X$ is algebraic. For instance, this functional
equation holds for any $X$ over a finite field. 

The main difference with Heinloth's formula
is that ours involves the \emph{determinant} of $M$, an invertible object of $\sM$ to be
defined in Definition \ref{d1}. This idea already appears in Andr\'e's book \cite[Proof
of 12.1.6.3]{andre}. For this, we use the theory of
$1$-dimensional objects developed by Kimura in \cite{kim2}. In the case of motives, we
conjecture that
$\det(M)$ is always the tensor product of a power of the Lefschetz motive and an Artin motive of
square $1$: we show that this conjecture is closely related to the conservativity of 
realisation functors on finite-dimensional motives, follows from the Hodge or the Tate
conjecture  (see Proposition \ref{p1}) and is true in important special cases (see Corollary
\ref{c2} and Remark \ref{r4} d)).

Theorem \ref{t1} almost provides a cycle-theoretic proof of the
functional equation for the usual zeta functions of motives over finite
fields. The catch is the sign conjecture. Unfortunately, the only known proof of this
conjecture over a finite field is by Katz-Messing \cite{km}, who rely on Deligne's proof of the
Riemann hypothesis \cite{weilII}\dots One could use Theorem \ref{t1} in the category of
$l$-adic sheaves, but this is not much more than the original approach of Grothendieck et al.

Over general fields, the situation is even more open. Nevertheless the class of smooth
projective varieties for which the sign conjecture is known is quite respectable: it
contains varieties of dimension $\le 2$, abelian varieties and complete intersections in a
projective space. It is also stable under products, as the full subcategory of $\sM(k)$
consisting of finite-dimensional objects is a (thick) tensor subcategory of
$\sM(k)$ by the results of \cite{kim1}. This puts the remark of
\cite[\S 6]{heinloth} in a broader context.

\subsection*{Acknowledgements} I wish to thank Noriyuki Otsubo for an enlightening discussion,
Yves Andr\'e for several comments on the first draft of this paper and Jean-Louis
Colliot-Th\'el\`ene for his help in the example of Remark \ref{r4} c).

\section{The determinant} Let $\sM$ be a $K$-linear tensor category, where $K$ is a field of
characteristic $0$. We assume that $\sM$ is
\emph{rigid} and that
$\End_\sM(\1)=K$. For any
$M$, we write
$\chi(M)=tr(1_M)\in K$; this is the \emph{Euler characteristic} of $M$. Unifying the
terminologies of Kimura and O'Sullivan, we shall say that $M$ is \emph{positive} (resp.
\emph{negative}) if there exists an integer $N>0$ such that $\Lambda^N(M)=0$ (resp.
$S^N(M)=0$), and \emph{finite-dimensional} if $M\simeq M^+\oplus M^-$ where $M^+$ is positive
and $M^-$ is negative. It can then be shown that in the positive (resp. negative) case,
$\chi(M)$ is a nonnegative (resp. nonpositive) integer and that the smallest such $N$ is
$\chi(M)+1$  (resp. $-\chi(M)+1$) (\cite{kim1}, \cite[9.1.7]{akos}). 

We denote by $\sM_\kim$ the full subcategory of finite-dimensional objects in $\sM$: it is
thick and rigid (\cite{kim1}, \cite[9.1.4 and 9.1.12]{akos}).

Recall that if $M$ is finite-dimensional, a decomposition $M\simeq M^+\oplus M^-$ is unique up
to isomorphism (\cite{kim1}, \cite[9.1.10]{akos}). This allows us to set
\begin{equation}\label{eq3}
\chi^+(M)=\chi(M^+),\qquad \chi^-(M) = \chi(M^-).
\end{equation}

\begin{defn}[\protect{\cite[15.2.3 and 15.2.4]{kim2}}]\label{d2} An object $L\in \sM$ is
\begin{itemize}
\item \emph{invertible} if there exists an object $T$ and an isomorphism $L\otimes T\simeq \1$;
\item \emph{$1$-dimensional} if it is either positive or negative, and $\chi(M)=1$ ($M$
positive) or $\chi(M)=-1$ ($M$ negative).
\end{itemize}
We denote by $\Pic(\sM)$ the full subcategory of $\sM$ consisting of invertible objects: it is
a symmetric monoidal subcategory of $\sM$.
\end{defn}

\begin{lemma}\label{l1} a) An object $L\in \sM$ is invertible if and only if it is
$1$-dimensional.\\ 
b) If $M$ is positive, $\Lambda^{\chi(M)}(M)$ is invertible.\\
c) If $M$ is negative, $S^{-\chi(M)}(M)$ is invertible.
\end{lemma}

\begin{proof} a) is the contents of \cite[15.2.6 and 15.2.9]{kim2}. In b) and c), it then
suffices to prove that these objects are $1$-dimensional. In case b), by \cite[7.2.4]{akos},
$\chi(\Lambda^n(M)) = \binom{\chi(M)}{n}$, and in case c) $\chi(S^n(M)) =
\binom{\chi(M)+n-1}{n}$ by the same reference.
\end{proof}

\begin{defn}[compare \protect{\cite[15.4.3]{kim2}}]\label{d1} Let $M\simeq M^+\oplus M^-$ be
finite dimensional. We define
\[\det(M) = \Lambda^{\chi^+(M)}(M^+)\otimes S^{-\chi^-(M)}(M^-)^{-1} \text{ (cf. \eqref{eq3}.)}\]
This is the \emph{determinant} of $M$: it is well-defined up to isomorphism.
\end{defn}

\begin{rks}\label{r3} a) If $L$ is invertible, we denote an object $T$ as in Definition \ref{d2}
by
$L^{-1}$. Note that $L^{-1}$ is unique up to unique isomorphism, and is the dual $L^*$ of
$L$.\\  b) $L\in \Pic(\sM)$ is positive (resp. negative) if and only if the switch of $L\otimes
L$ is trivial (resp. equals $-1$).\\ c) If $M$ is positive, $\det(M)$ is positive. If $M$ is
negative, $\det(M)$ has the sign $(-1)^{\chi(M)}$: this follows from the proof of Lemma
\ref{l1} c).
\end{rks}

The following lemma will be used in the proof of Proposition \ref{p1}.

\begin{lemma}\label{l5} Let $f:M_1\to M_2$ be a morphism between finite-dimen\-sion\-al objects
such that $\chi^+(M_1) = \chi^+(M_2)$ and $\chi^-(M_1) = \chi^-(M_2)$. Then $f$ induces a
morphism $\det(f):\det(M_2)\to \det(M_2)$. (In view of Definition \ref{d1}, we don't claim that
$\det(f)$ is canonically defined unless $M$ is purely positive or negative). Moreover, $f$ is
an isomorphism if and only if $\det(f)$ is.
\end{lemma}

\begin{proof} Choose decompositions $M_i =M_i^+\oplus M_i^-$, and let
\[f = \begin{pmatrix}
f^{++}& f^{-+}\\
f^{+-}& f^{--}
\end{pmatrix}\]
be the corresponding matrix description of $f$. We set $\det(f)=\det(f^{++})\allowbreak\otimes \det(f^{--})^{-1}$. If $f$ is an isomorphism, then $f^{++}$ and $f^{--}$ are isomorphisms. To see this, one may go modulo the ideal $\sN$  of morphisms universally of trace $0$: by \cite[9.1.14]{akos}, $\sN$ is locally nilpotent and $f^{+-},f^{-+}\in \sN(M_1,M_2)$ by the proof of this proposition. 

For the converse, we work in $\sM_\kim$. By \cite[9.2.2]{akos},  $\sM_\kim/\sN$ is abelian semi-simple. Let $\bar f$ be the image of $f$ in this category: we may write $\bar f$ as the direct sum of an isomorphism $f':M'_1\to M'_2$ and a $0$ morphism $f'':M''_1\to M''_2$. Then $\det(f) = \det(f')\otimes \det(f'')$. But if $M''_1$ and $M''_2$ were nonzero, we would have $\det(f'')=0$, hence $\det(\bar f)=0$, contrary to the hypothesis. Thus $\bar f$ is an isomorphism, and by nilpotency $f$ is an isomorphism.
\end{proof}

\begin{lemma} \label{l2} a) If $M$ is positive, we have isomorphisms
\[\Lambda^n(M^*)\simeq \Lambda^{\chi(M)-n}(M)\otimes \det(M)^*\]
for all $n\in [0,\chi(M)]$.\\
b) If $M$ is negative, we have isomorphisms
\[S^n(M^*)\simeq S^{-\chi(M)-n}(M)\otimes \det(M)^*\]
for all $n\in [0,-\chi(M)]$.
\end{lemma}

\begin{proof} a)  Note that $\Lambda^n(M^*) \simeq \Lambda^n(M)^*$. The obvious pairing
\[\Lambda^n(M)\otimes \Lambda^{\chi(M)-n}(M)\to \det(M)\]
and its dual applied with $M^*$ instead of $M$ define morphisms
\begin{align*}
\Lambda^n(M)\otimes \left(\Lambda^{\chi(M)-n}(M)\otimes \det(M)^*\right)\to \1\\
\1\to \left(\Lambda^{\chi(M)-n}(M)\otimes \det(M)^*\right)\otimes\Lambda^n(M)
\end{align*}
and it is easy to check that they verify the duality axioms \cite{dp}.
The proof for b) is similar.
\end{proof}

\section{The functional equation} We can now state our theorem:

\begin{thm}\label{t1} Let $M\simeq M^+\oplus M^-$ be finite dimensional, and $M^*$ be its dual.
Then we have:
\[Z_\sM(M^*,T^{-1})=(-1)^{\chi^+(M)} \det(M) T^{\chi(M)} Z_\sM(M,T).\]
\end{thm}

\begin{proof} Since the zeta function transforms direct sums into products, we may treat
separately the cases where $M$ is positive and negative. Note that $M^*$ has the same sign as
$M$.

1. If $M$ is positive, we have (see \cite[proof of 13.3.3.1]{andre}):
\[Z_\sM(M^*,T^{-1})^{-1}= \sum_{n= 0}^{\chi(M)} [\Lambda^n(M^*)] (-T)^{-n}.\]

Applying Lemma \ref{l2} a), we may rewrite this sum as
\begin{multline*}
\relax  [\det(M)^*]\sum_{n= 0}^{\chi(M)} [\Lambda^{\chi(M)-n}(M)] (-T)^{-n}\\
= [\det(M)^*]\sum_{n= 0}^{\chi(M)} [\Lambda^{n}(M)] (-T)^{n-\chi(M)}\\
= [\det(M)^*](-T)^{-\chi(M)}\sum_{n= 0}^{\chi(M)} [\Lambda^{n}(M)] (-T)^{n}\\
= [\det(M)^*](-T)^{-\chi(M)}Z_\sM(M,T)^{-1}.
\end{multline*}

2. If $M$ is negative, we have
\[Z_\sM(M^*,T^{-1}) =\sum_{n= 0}^{-\chi(M)} S^n(M^*) T^{-n}.\]

Applying Lemma \ref{l2} b), we may rewrite this sum as
\begin{multline*}
\relax  [\det(M)^*]\sum_{n= 0}^{-\chi(M)} [S^{-\chi(M)-n}(M)] T^{-n}\\
= [\det(M)^*]\sum_{n= 0}^{-\chi(M)} [S^{n}(M)] T^{n+\chi(M)}\\
= [\det(M)^*]T^{\chi(M)}\sum_{n\ge 0} S^{n}(M)] T^{n}\\
= [\det(M)^*]T^{\chi(M)}Z_\sM(M,T).
\end{multline*}
\end{proof}

\begin{rks}\label{r2} a) As remarked in \cite[13.2.1.1]{andre}, the map $K_0(\sM_\kim)\to
K_0(\sM_\kim/\sN)$ is an isomorphism thanks to the local nilpotency of $\sN$: therefore it
makes no difference to work in $\sM_\kim$ or in $\sM_\kim/\sN$ as long as the zeta function is
concerned.\\ b) One may wonder if $\chi^-(M)$ is always even: this would give a nicer formula
in Theorem \ref{t1}. As Yves Andr\'e pointed out, this is of course false if we take for $\sM$
the category of $\Z/2$-graded $\Q$-vector spaces. However, this turns out to be true (most of
the time, and conjecturally always) if $\sM$ is the category of motives with rational
coefficients over a field: Yves Andr\'e gave us an argument for this in characteristic $0$, and
his argument inspired part of the proofs given in the next section. 
\end{rks}

\section{The case of motives} We now assume that $\sM=\sM_\sim(k)$ is the category of motives
over a field $k$ with respect to a given adequate equivalence relation $\sim$: we shall often
abbreviate $\sM_\sim(k)$ into $\sM(k)$.  The two main examples for $\sim$ are rational and
homological equivalence. We only work with finite-dimensional motives $M$ relatively to $\sim$.
If $\sim$ is homological equivalence, this hypothesis is equivalent to the sign conjecture. 

Notice in any case that, in view of Remark \ref{r2} a), the motivic zeta function of $M$ with
respect to $\sim$ is the same as its motivic zeta function with respect to any coarser adequate
equivalence relation. Similarly, if $H$ is a classical Weil cohomology theory and $\sim$ is
coarser than $H$-equivalence, then $K_0(\sM_H(k))\iso K_0(\sM_\sim(k))\iso K_0(\sM_\num(k))$ by
the nilpotence result of \cite[Prop. 5]{ak}, without even restricting to finite-dimensional
motives. Another such example is Voevodsky's nilpotence theorem \cite{V}, which implies that
$K_0(\sM_\rat(k))\iso K_0(\sM_\alg(k))$.

Let $L\in \sM_\sim(k)$ be an invertible motive. Let $L_\num$ be the image of $L$ in
$\sM_\num(k)$. By \cite[Prop. 5]{ak}, we may lift $L_\num$ to an invertible motive
$L_H\in\sM_H(k)$ with respect to any classical Weil cohomology theory, uniquely up to
isomorphism. 

\begin{defn} The \emph{weight} of $L$ is the weight of the $H$-realisation of $L_H$ under a
suitable classical Weil cohomology theory $H$.
\end{defn}

In characteristic $0$, a ``suitable" Weil cohomology theory is either Hodge cohomology or
$l$-adic cohomology; in characteristic $p$, it is either crystalline cohomology or $l$-adic
cohomology. That the weight does not depend on the choice of the suitable cohomology theory
follows from Artin's comparison theorem in characteristic $0$ and from Deligne's proof of the
Riemann hypothesis in characteristic $p$ \cite{km}.

\begin{conj}\label{c3} Let $k$ be a field. There is no invertible motive $L$ of odd weight $n$
in $\sM(k)$.
\end{conj}

Equivalently: there is no negative invertible motive in $\sM(k)$.

In fact, we can almost prove Conjecture \ref{c3}:

\begin{prop}\label{l4} Conjecture \ref{c3} is true in the following cases:
\begin{enumerate}
\item $\car k=0$.
\item $k$ is a regular extension of $\F_q$, where $q$ is not a square.
\item $n=1$ and $L$ is effective.
\end{enumerate}
Moreover, the Tate conjecture over finite fields implies Conjecture \ref{c3} in the remaining
cases. (Can one avoid the recourse to the Tate conjecture???)
\end{prop}

\begin{proof} Clearly, we may work modulo homological equivalence.

(1) This is obvious via the Hodge realisation: in fact, any
invertible pure Hodge structure is of the form $\Q(r)$ since the Hodge numbers go by symmetric
pairs. (This was part of Andr\'e's argument alluded to in Remark \ref{r2}.)  Alternately, we may
use
$l$-adic cohomology as follows: Let $k_0\subseteq k$ be a finitely generated field over
which $L$ is defined. Let $S$ be a regular model of $k_0$ of finite type over $\Z$. Then one may
find a closed point $s\in S$ with residue field $\F_p$ ($p$ a large enough prime number) such
that $L$ has good reduction at $s$. This reduces us to the case where
$k=\F_p$. Let $F$ be the Frobenius endomorphism of $L$. Since $\End(\1)=\Q$, $tr(F)$ is a
rational number, which is also the (unique) eigenvalue of the action of $F$ on $H_l(L)$.
But if $n$ is the weight of $L$, $H_l(L)$ is a direct summand of $H^n(\bar X,\Q_l)(r)$ for some
smooth projective variety $X$ and some $r\in \Z$; by the Weil conjecture, $|tr(F)| = p^{n/2-r}$,
a contradiction.

(2) The same $l$-adic argument as in (1) works, and a similar argument works for crystalline
cohomology.

(3) Write $L$ as a direct summand of $h(X)$ for some smooth projective variety $X$. We may write
$L$ as a direct summand of $h^1(X)$. (Note that there is always a canonical decomposition
$h(X)=h^0(X)\oplus h^1(X)\oplus h^{\ge 2}(X)$ and that $\Hom(L,h^0(X))=\Hom(L,h^{\ge 2}(X))=0$
for weight reasons.) Replacing $X$ by its Albanese variety, we may assume that $X$ is an abelian
variety $A$. Now, by Poincar\'e's complete irreducibility theorem, the full subcategory of
motives modulo homological equivalence consisting of the $h^1(A)$ is abelian semi-simple; thus
$L=h^1(B)$ for some abelian variety $B$. But this is impossible, since $L$ is $1$-dimensional
and $\dim h^1(B)$ has to be even.

Finally, if the Tate conjecture holds over $\F_{q^2}$, then the generalised Tate conjecture also
holds thanks to Honda's theorem (compare \cite[proof of Th. 2]{kFq}); hence, up to twisting $L$ by some power of the Lefschetz
motive, we may assume that it is effective of weight $1$ and we are reduced to (3).
\end{proof}

\begin{rk} Notice that Conjecture \ref{c3} becomes
false over $\F_{q^2}$ after extending scalars from $\Q$ to a suitable quadratic extension, as
shown by $h^1$ of a supersingular elliptic curve. Similarly, Conjecture \ref{c3} becomes (very)
false over any field after extending scalars from $\Q$ to $\bar \Q$.
\end{rk}

\begin{thm}\label{t2} a) If conjecture \ref{c3} is true for $k$ (compare Proposition \ref{l4}),
$\chi(M)$ is even for any negative $M\in \sM(k)$. Thus, in this case, the functional equation
of Theorem
\ref{t1} simplifies as
\[Z_\sM(M^*,T^{-1})=\det(M) (-T)^{\chi(M)} Z_\sM(M,T).\]
b) $\chi^-(M)$ is even over any field $k$ if $M=h(X)$ for $X$ a smooth projective variety
verifying the sign conjecture. We then have
\begin{equation}\label{eq1}
Z_\sM(h(X),T^{-1}) = (-T)^{\chi(X)} \det(h(X))Z_\sM(h(X), [\L]^{-d}T).
\end{equation}
This formula holds modulo homological equivalence, and even modulo rational equivalence if $h_\rat(X)$ is finite-dimensional.
\end{thm}

\begin{proof} a) If there exists a negative motive $M$ of odd Euler characteristic, then
$\det(M)$ is negative of Euler characteristic $-1$ (see Remark \ref{r3} c)).
It remains to justify b): this is true since the odd-dimensional $l$-adic
Betti numbers of $X$ are even by Hard Lefschetz \cite[(4.1.5)]{weilII}. The functional equation
for $h(X)$ then follows from the formula $h(X)^*=h(X)\otimes
\L^{-d}$, where $\L$ is the Lefschetz motive, and the obvious identity
\[Z_\sM(h(X)\otimes \L^{-d},T)= Z_\sM(h(X), [\L]^{-d}T).\]
\end{proof}

\section{The determinant of a motive}

When $X$ is an abelian variety, \eqref{eq1} is Heinloth's formula, with the difference that the
factor $\det(h(X))$ is replaced by $1$. We want to explain in this section why this happens. This is a good opportunity to introduce the following
strengthening of Conjecture \ref{c3}:

\begin{conj}\label{c1} For any field $k$, $\Pic(\sM_\sim(k))$ is generated by $\L$ and
Artin motives of square $1$.
\end{conj}

(Such Artin motives are in $1$-to-$1$ correspondence with the elements of $k^*/k^{*2}$.)

Recall that a functor $F$ is \emph{conservative} if, when $F(f)$ is an isomorphism, $f$ is
already an isomorphism. We shall say that $F$ is \emph{essentially injective} if it is
injective on isomorphism classes of objects. If $F$ is fully faithful, it is both
conservative and essentially injective.

\begin{prop}\label{p1} a) Suppose that $\sim$ is finer than homological equivalence with
respect to some Weil cohomology theory $H$. Then conjecture \ref{c1} implies that the
realisation functor $H$ is conservative on finite-dimen\-sion\-al motives.\\  
b) If $k$ verifies Conjecture \ref{c3}, the essential injectivity of the Hodge realisation (over
a field of characteristic $0$) or the $l$-adic realisation (over a finitely generated field)
implies Conjecture \ref{c1} for any $\sim$. In particular, Conjecture \ref{c1} follows from
either the Hodge or the Tate conjecture.
\end{prop}

\begin{proof} a) Let $f:M_1\to M_2$ be a morphism between finite-dimensional motives such that
$H(f)$ is an isomorphism for some Weil cohomology theory $H$. This implies that
$\chi^\pm(M_1)=\chi^\pm(M_2)$, hence, by Lemma \ref{l5}, $f$ induces a morphism
$\det(f):\det(M)\to \det(M')$. Obviously, $H(\det(f))=\det H(f)$, hence $H(\det(f))$ is an
isomorphism. But under Conjecture \ref{c1}, $H$ is clearly conservative on invertible motives,
hence $\det(f)$ is an isomorphism and $f$ is an isomorphism by Lemma \ref{l5}.

For b), we may reduce to the case where $\sim$ is finer than Hodge or $l$-adic equivalence,
thanks to \cite[Prop. 5]{ak}.

1) Let us assume that $\car k=0$ and that the Hodge realisation is conservative. We may assume
that $L\in \Pic(\sM(k))$ is defined over a subfield $k_0$ of
$\C$. As seen in the proof of Proposition
\ref{l4},
$H(L)\simeq H(\L^n)$ for some
$n\in\Z$; by essential injectivity, $L_\C\simeq \L^n$, hence up to twisting we may assume that
$L_\C=\1$. There is a finitely generated extension $k_1$ of $k_0$ such that $L_{k_1}=\1$, and
after a suitable specialisation we may assume that $k_1/k_0$ is finite. Taking traces, we get
that $L\otimes h(k_1)\simeq h(k_1)$ in $\sM(k_0)$; hence $L$ is a direct summand of $h(k_1)$
and is an Artin motive. But Artin motives are in $1$-to-$1$ correspondence with rational Galois
representations of $k$, and any rational character of a (pro)finite group has square $1$.

2) Let us assume that $k$ is finitely generated, verifies Conjecture \ref{c3} and that the
$l$-adic realisation is conservative. Let $L$ be an invertible motive and let $R_l(L)$ be its
$l$-adic realisation. Let $n$ be the weight of $L$, which is assumed to be even. Using the same
rationality argument as in the proof of Proposition \ref{l4} (1), we see that the eigenvalue of
Frobenius acting on the specialisation of $R_l(L)$ at any finite place of good reduction and
residue field $\F_q$ is of the form $\pm q^{n/2}$. Thus Frobenius automorphisms act trivially
on the $l$-adic representation $R_l(L^{\otimes 2}\otimes \L^{-n})$; by \v Cebotarev's density
theorem, it follows that $R_l(L^{\otimes 2}\otimes \L^{-n})$ is trivial, hence $R_l(L\otimes
\L^{-n/2})$ is given by a quadratic character of $k$. Let $A$ be the corresponding Artin motive
of square $1$: by essential injectivity we have $L\simeq A\otimes L^{n/2}$, as desired. 

It remains to justify the last assertion of b). If the Hodge conjecture holds, then the Hodge
realisation functor is essentially injective on the level of motives modulo homological
equivalence. But let $L,L'\in \Pic(\sM_\sim(k))$: since $L$ and $L'$ are finite-dimensional,
Kimura's nilpotence theorem implies that any isomorphism between them modulo homological
equivalence lifts to an isomorphism modulo $\sim$. The argument is the same with the Tate
conjecture.
\end{proof}

\begin{rks} 1) With a little more effort, one can see that conjecture \ref{c1} implies the
irreducibility of the motive $\1$ modulo homological equivalence. Indeed, let $f:L\to \1$ be a
monomorphism: then
$H(f)$ is also a monomorphism, hence $L=0$ or $\chi(L) = \pm 1$ and one concludes immediately.
The conservativity in Proposition \ref{p1} a) then follows from 
\cite[12.1.6.3]{andre}: I thank Yves Andr\'e for pointing this out. In fact, after reading the
proof \cite[12.1.6.3]{andre}, I realised that it is very closely related to the proof of Prop.
\ref{p1} a) given above!\\
2) Conjecture \ref{c1} is obviously equivalent to the following: for any finite-dimensional
$M\in\sM(k)$, we have
\begin{equation}\label{eq2}
\det(M)=\L^r\otimes A
\end{equation} 
for some $r\in\Z$, where $A$ is a $1$-dimensional Artin motive of square $1$.
It is interesting to investigate this reformulation in the light of formal
propreties of $\det$, as follows:
\end{rks}

\begin{prop}\label{l3} Let $\sM$ be a rigid tensor $\Q$-linear category. Then, for any
finite-dimensional $M,N$:
\begin{enumerate}
\item $\det(M\oplus N) =\det(M)\det(N)$.
\item $\det(M\otimes N) = \det(M)^{\chi(N)}\otimes \det(N)^{\chi(M)}$.
\item $\det(\Lambda^n(M))= \det(M)^{r}$, with $r=n\binom{\chi(M)}{n}/\chi(M)$.
\item $\det(S^n(M)) = \det(M)^{s}$, with $s=n\binom{\chi(M)+n-1}{n}/\chi(M)$.
\end{enumerate}
\end{prop}

\begin{proof} This follows readily from identities between Schur functors as in
\cite[\S 4]{heinloth}.
\end{proof}

\begin{cor}\label{c2} a) the subset $K'_0(\sM(k))$ of $K_0(\sM(k))$ consisting of differences
of classes of finite-dimensional motives that verify \eqref{eq2} is a
sub-$\lambda$-ring of $K_0(\sM(k))$.\\ 
b) Similarly, the subset $K''_0(\sM(k))$ of $K_0(\sM(k))$ consisting of differences of classes of finite-dimensional motives such that $A=1$ in \eqref{eq2} is a sub-$\lambda$-ring of $K_0(\sM(k))$.\\ 
c) Suppose that $M$ is weakly polarisable in the sense that $M\simeq
M^*\otimes \L^n$ for some integer $n$. Then $\det(M)^2=
\L^{n\chi(M)}$.\\
d) Suppose that $X/k$ verifies the standard conjecture B. Then, for all $i\ge 0$, $\det(h^i(X))^2=\L^{(-1)^iib_i(X)}$, where $b_i(X)$ is the $i$-th Betti number of $X$. Moreover, letting $d=\dim X$:
\[\det(h(X))=
\begin{cases}
\L^{d\chi(X)/2}&\text{if $d$ is odd}\\
\L^{\frac{d}{2}(\chi(X)-(-1)^{d/2}db_{d/2}(X))}\otimes \det(h^{d/2}(X))&\text{if $d$ is even.}
\end{cases}
\]
e) $K''_0(\sM(k))$ contains $[h^1(X)]$ for any smooth projective variety $X$. (Recall that the
direct summand $h^1(X)$ of the Chow motive $h(X)$ was constructed by Murre \cite{Mu1}).\\
f) If $C$ is a curve, $\det(h(C))= \L^{\chi(C)/2}$.\\
g) If $A$ is an abelian variety, $\det(h(A))=1$.
\end{cor}

\begin{proof} a), b) and c) follow directly from Proposition \ref{l3} (see also Remark
\ref{r3} a) for c)). The first claim of d) follows from c) (note that
$\chi(h^i(X))=(-1)^ib_i(X)$). We get the second one by grouping the terms $\det(h^i(X))$ and
$\det(h^{2d-i}(X))$ together for $i<d/2$, and using that $h^{2d-i}(X)\simeq h^i(X)^*\otimes
\L^{d-i}$.

In e), note that $h^1(X)\simeq h^1(A)$, where $A$ is the Albanese variety of $X$; the claim
then follows from the results of Shermenev, Beauville and Deninger-Murre, saying that
$h^i(A)\simeq S^i(h^1(A))$ (apply this for $i=2\dim A$). f) follows from e). Finally, in g),
$A$ verifies the standard conjecture B by Lieberman-Kleiman; if $\dim A = g$, we have $b_i(A) =
\binom{g}{i}$ and the claim follows from b) and d).
\end{proof}

\begin{rks}\label{r4}
a) Theorem \ref{t2} and Corollary \ref{c2} f) give back Kapranov's functional equation for
curves \cite{kapranov},  with coefficients in $K_0(\sM(k))$ (a little information is lost). \\
b) The reader familiar with Heinloth's paper will recognize part of her arguments in the proof
of Corollary \ref{c2} g). This corollary was actually inspired by reading her paper.\\ 
c) In view of Corollary \ref{c2} f) and g), one might expect that the Artin motive $A$ appearing
in \eqref{eq2} is always $\1$ for $M$ of the form $h(X)$. This is wrong: rational surfaces give
examples where  \eqref{eq2} holds but $A$  is nontrivial. Indeed, for such a surface $S$,
$h^2(S)= \L\otimes NS_S$, where $NS_S$ is the Artin motive corresponding to the N\'eron-Severi
lattice of $S$. Hence $\det(h^2(S))= \L^\rho \otimes\det (NS_S)$, where $\rho$ is the Picard
number. For examples where the action of the absolute Galois group of $k$ on $NS_S$ has
nontrivial determinant $d\in k^*/k^{*2}$, we may take for $S$ the blow-up of $\P^2$ at
$\{\sqrt{d}, -\sqrt{d}\}$ (I thank Colliot-Th\'el\`ene for pointing this out): this shows that
all Artin motives of square $1$ are caught thusly. This also incidentally explains the sign
$(-1)^\mu$ in \cite[p. 9]{kl}. I expect that this is typical of how nontrivial Artin motives
may occur.\\ 
d) The class of varieties $X$ such that $h(X)\in K'_0(\sM(k))$ is closed under
products; it is also closed under projective bundles and blow-ups with smooth centres, in the
sense that of $[h(X)]\in K'_0(\sM(k))$ and
$Z$ is a closed smooth subvariety of $X$ such that $[h(Z)]\in K'_0(\sM(k))$, then
$[h(Bl_Z(X))]\in K'_0(\sM(k))$. Thus surfaces $S$ such that $h(S)$ verifies \eqref{eq2}
include rational surfaces, ruled surfaces, products of two curves and abelian surfaces. The
first open case for Conjecture \ref{c1} seems to be that of a K3 surface.
\end{rks}

\end{document}